\documentclass[11pt]{article}
\usepackage{amssymb,latexsym,amsmath} 
\usepackage{dsfont}
\usepackage{upgreek}
\usepackage{graphics}
\usepackage{graphicx}
\usepackage{lscape}
\usepackage{multirow}
\usepackage{amssymb}
\usepackage{dsfont}
\usepackage{authblk}
\usepackage[switch,columnwise]{lineno}
\usepackage{lmodern}
\usepackage{fourier}
\usepackage{float,lscape}
\usepackage{multirow}
\usepackage{array}
\usepackage{book tabs}
\usepackage{subcaption}
\usepackage[labelformat=parens,labelsep=quad,skip=3pt]{caption}
\usepackage{graphicx}
\usepackage[margin=1in]{geometry}
\usepackage{booktabs}
\date{}

\begin{document}

\title{Exponentiated Inverse Power Lindley Distribution and its Applications}
\author[1*]{Rameesa Jan}
\author[2]{T.R.Jan}
\author[3]{ Peer Bilal Ahmad}
\affil[1] {Department of Statistics,  University Kashmir, bhatrumaisa2@gmail.com,*\underline{corresponding author}}
\affil[2] {Department of Statistics,  University of Kashmir, drtrjan@gmail.com}
\affil[3] {Department of Mathematical Sciences,IUST, peerbilal@yahoo.co.in}
\maketitle
\paragraph{Abstract} \mbox{}\\
In this article,a three parameter generalization of inverse lindley distribution is obtained, with the purpose of obtaining a more flxeible model relative to the behaviour of hazard rate functions. Various statistical properties such as density, hazard rate functions, moments, moment generating functions, stochastic ordering, renyi entropy, distribution of rth order statistics has been derived. The method of maximum likelihood estimation has been used to estimate parameters. Further confidence intervals are also obtained. Finally applicability of the proposed model to the real data is analyzed. A comparison has also been made with some existing distributions.	

	\section{Introduction}
One of the most important branch of Statistics is survival and reliability analysis. Survival analysis deals with analyzing the expected duration of time until one or more event happens such as death in biological organisms and failure in mechanical systems. Survival analysis has its applications in various applied sciences such as engineering, actuarial science, demography, public health and industrial reliability and is usually used to describe the survival times of members of a group through life tables, survival functions and hazard rates etc. For modeling survival data various models has been defined based on the behavior of hazard rate (monotone or non monotone).\\
Most real-life systems exhibit non-monotone (bathtub and upside down bathtub) shape for their hazard rates. For analysis of lifetime data having bathtub-shaped hazard rate function various probability distributions have been proposed. Mudholkar and Srivastav $(1993)$, Xie and Lai $(1996)$, Xie and T.Tang $(2002)$ introduced several extensions of weibull distribution to model bathtub-shaped hazard rate date.The inverse versions of the probability distributions are able to model the data that show upside-down bathtub-shape for their hazarad rates e.g inverse weibull distribution, inverse gamma distribution,inverse guassian distribution.\\
   Lindley $(1958)$ introduced a one parameter distribution known as Lindley distribution with probability density function(pdf) and cumulative distribution function(CDF) respectively given by:  
\begin{eqnarray}
f(x,\beta)&=&\frac{\beta^2}{1+\beta}\left(1+x\right)e^{-\beta x} \; ;  \; x>0,\beta>0 \nonumber \\
F(x,\beta)&=&1-\left(1+\frac{\beta x}{1+\beta}\right)e^{-\beta x} \nonumber
\end{eqnarray}
Ghitany, Atich, Nadarajah $(2008)$ studied various statistical properties of Lindley dsitribution and showed the superiority of lindley distribution over exponential distribution for waiting times before service of bank customers. The lindley distribution has been extended by different researchers including Zakerzadeh and Dolati$(2009)$,Nadarajah S, Bakouch HS, Tahmasbi R $(2011)$, Bakouch HS, Al-Zaharani B, Al-Shomrani A, Marchi V, Louzad F $(2012)$, Shanker and Mishra$(2013)$, Ghitany M, Al-Mutairi D, Balakrishnan N, Al-Enezi I $(2013)$, Ashour and Eltehiwy$(2015)$.
Sharma,Singh,Singh,Agiwal $(2015)$ introduced the inverted version of lindley distribution known as Inverse lindley distribution with pdf and cdf respectively given by:\\
\begin{eqnarray}
f(x,\beta)&=&\frac{\beta^2}{1+\beta}\left(\frac{1+x}{x^3}\right)e^{\frac{-\beta}{x}} \; ; \;  x>0,\beta>0  \\
F(x,\beta)&=&\left(1+\frac{\beta}{1+\beta}\frac{1}{x}\right)e^{\frac{-\beta}{x}} \nonumber
\end{eqnarray}
Sharma V, Singh S, Singh U, Merovci F $(2015)$ added another parameter to the inverse lindley distributions by taking the inverse power transformation of a random variable with lindley distribution and named it as Generalized inverse lindley distribution. Note that Barco, Mazucheli, Janerio also obtained the generalised inverse lindley distribution by taking the transformation $X=Y^{\frac{1}{\alpha}}$ where $Y$ follows inverse lindley distribution known as Inverse power lindley distribution with pdf given as in equation (2).\\   
\begin{eqnarray}
f(x,\beta)&=&\frac{\alpha \beta ^2}{1+\beta}\left(\frac{1+x^\alpha}{x^{2\alpha+1}}\right)e^{\frac{-\beta}{x^{\alpha}}}  \;  ;  \;x>0 ,\alpha>0  ,\beta>0\\
F(x,\beta)&=&\left(1+\frac{\beta}{1+\beta}\frac{1}{x^\alpha}\right)e^\frac{-\beta}{x^\alpha}  \nonumber
\end{eqnarray}\\

A new three parameter probability distribution using the transformation $H(z)=[G(z)]^{\theta}$, where $G(z)$ is a CDF and $\theta$ is a positive real number, has been proposed. The new distribution Exponentiated Inverse power lindley distribution (EIPLD) thus obtained involves inverse power lindley distribution and inverse lindley distribution as its sub-models. Some mathematical properties of the proposed distribution has been derived.The moment generating function, quantile function, order statistics, renyi entropy and stochastic ordering has also been discussed.Maximum likelihood estimation of parameters and their confidence intervals are derived.Finally the proposed distribution is analyzed using the real life data and is compared to other well-known existing distributions.\\
\;

The proposed distribution is most conveniently specified in terms of cumulative distribution function given as:
\begin{equation}
G_{\theta}(z)=\left[\left(1+\frac{\beta}{1+\beta}\frac{1}{z^\alpha}\right)e^{\frac{-\beta}{z^{\alpha}}}\right]^{\theta} 
\end{equation}

and the corresponding probability distribution function is given by

\begin{equation}
g_{\theta}(z)=\frac{\alpha\beta^{2}\theta}{1+\beta}\left(\frac{1+z^\alpha}{z^{2\alpha+1}}\right)e^{\frac{-\theta\beta}{z^\alpha}}\left[1+\frac{\beta}{1+\beta}\frac{1}{z^{\alpha}}\right]^{\theta-1}      
\end{equation} 
\rlap{\hspace*{35em};$z>0$ , $\alpha >0$ ,$\beta >0$, $\theta >0$}\\
\textbf{Special Cases}:\\
Case 1: For $\theta=1$ EIPLD reduces to inverse power lindley distribution.\\
Case 2: For $\alpha=\theta=1$ EIPLD reduces to inverse lindley distribution.

\section{Various Reliability Measures}

The survival function of EIPLD is given by (5) 
\begin{equation}
S_{z}(z)=1-\left[\left(1+\frac{\beta}{1+\beta}\frac{1}{z^\alpha}\right)e^{\frac{-\beta}{z^{\alpha}}}\right]^{\theta} 
\end{equation}\\

and the corresponding hazard function of EIPLD can be obtained as $h_{z}(z)=\frac{g_{z}(z)}{S_{z}(z)}$ and is given by (6)
\begin{equation}
h_{Z}(z)=\frac{\frac{\alpha\beta^{2}\theta}{1+\beta}\left(\frac{1+z^\alpha}{z^{2\alpha+1}}\right)e^{\frac{-\theta\beta}{z^\alpha}}\left[1+\frac{\beta}{1+\beta}\frac{1}{z^{\alpha}}\right]^{\theta-1} }{1-\left[\left(1+\frac{\beta}{1+\beta}\frac{1}{z^\alpha}\right)e^{\frac{-\beta}{z^{\alpha}}}\right]^{\theta}}  
\end{equation}\\

and the proportional reversed hazard rate function of EIPLD is given by
\begin{eqnarray}
\lambda(z)&=&\frac{d}{dz}\ln G_{\theta}(z)=\frac{g_{\theta}(z)}{G_{\theta}(z)}\nonumber \\
          &=&\frac{\frac{\alpha\beta^{2}\theta}{1+\beta}\left(\frac{1+z^\alpha}{z^{2\alpha+1}}\right)e^{\frac{-\theta\beta}{z^\alpha}}\left[1+\frac{\beta}{1+\beta}\frac{1}{z^{\alpha}}\right]^{\theta-1} }{\left[\left(1+\frac{\beta}{1+\beta}\frac{1}{z^\alpha}\right)e^{\frac{-\beta}{z^{\alpha}}}\right]^{\theta}}  
\end{eqnarray}
 
The pdf plot, survival function,hazard rate and proportional hazard rate for different values of parameters are presented in Fig 1,Fig 2,Fig 3,Fig 4 respectively.

	\begin{figure}[h]
		\begin{subfigure}{6cm}
			\centering\includegraphics[width=6cm]{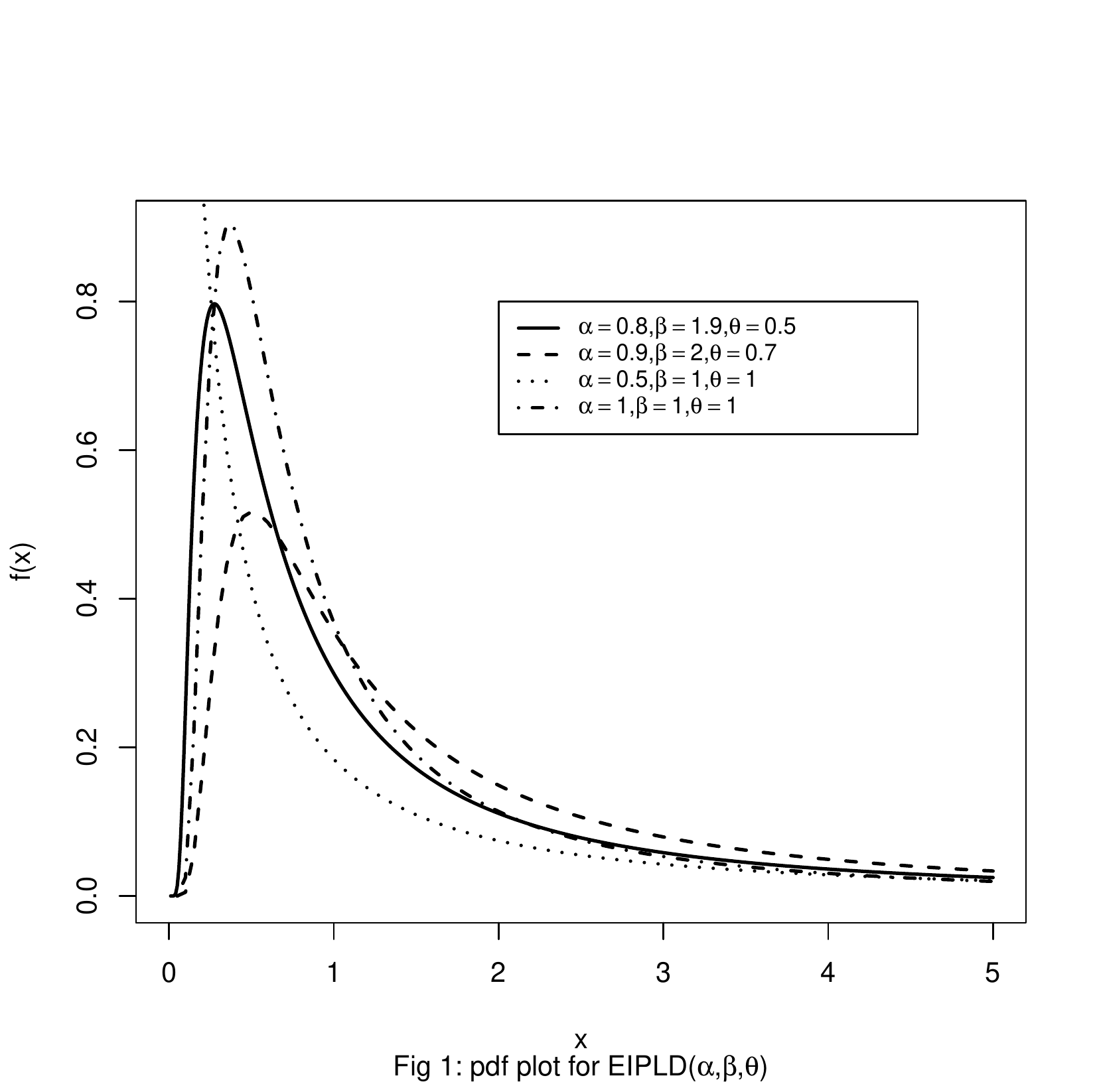}
		
		\end{subfigure}
		\begin{subfigure}{6cm}
			\centering\includegraphics[width=6cm]{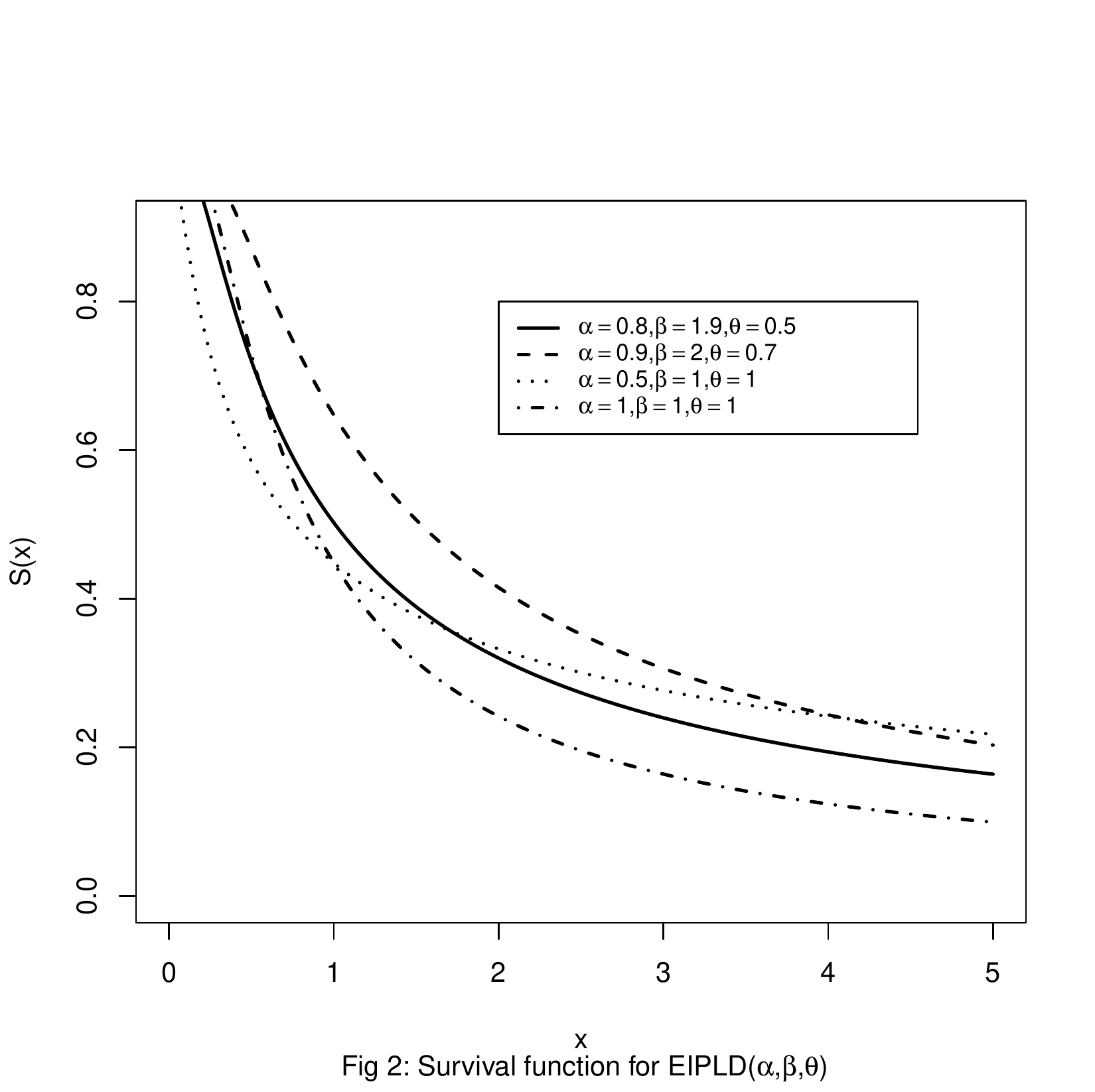}
		
		\end{subfigure}
		
		\begin{subfigure}{6cm}
			\centering\includegraphics[width=6cm]{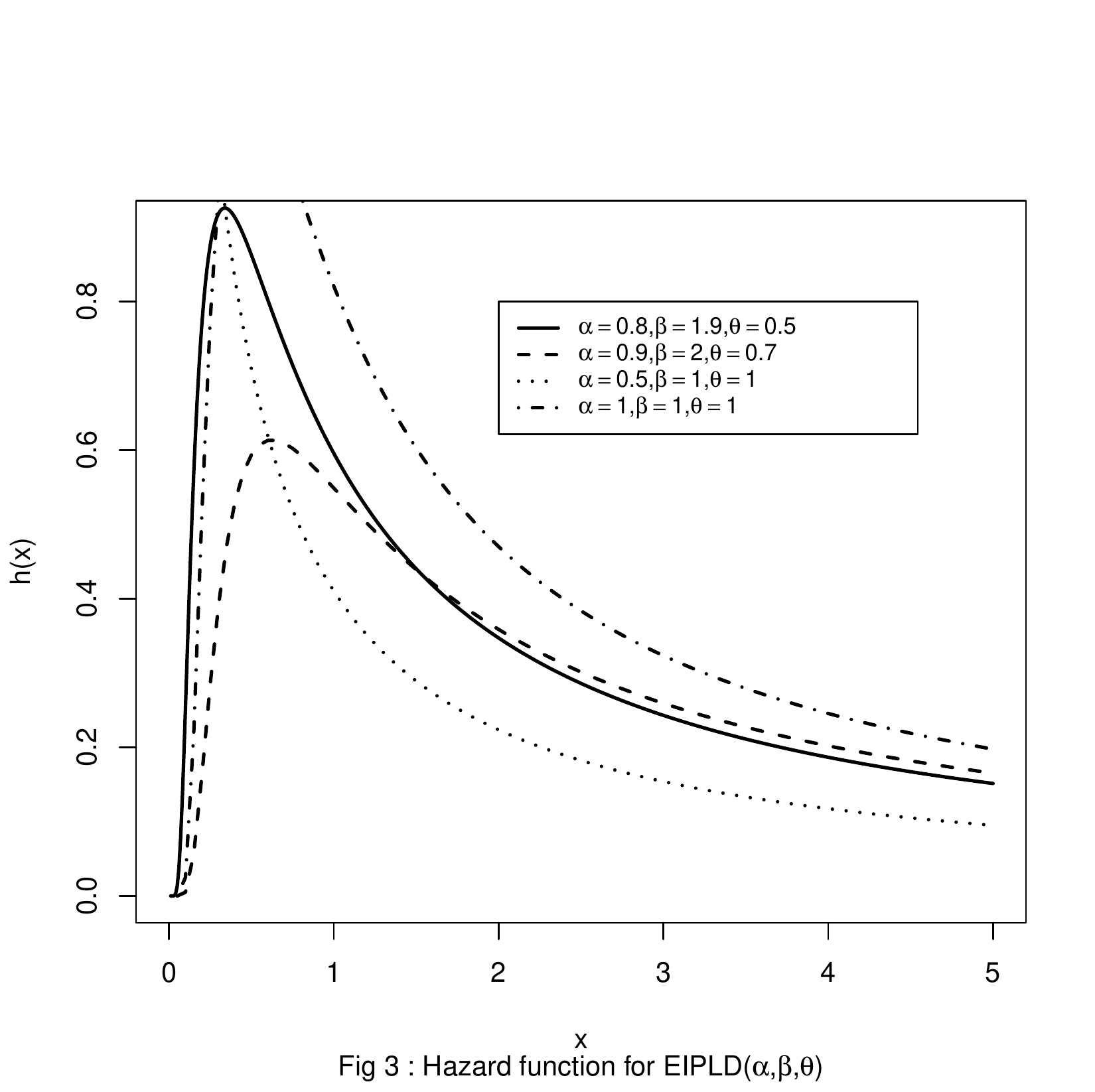}
		
		\end{subfigure}
		\begin{subfigure}{6cm}
			\centering\includegraphics[width=6cm]{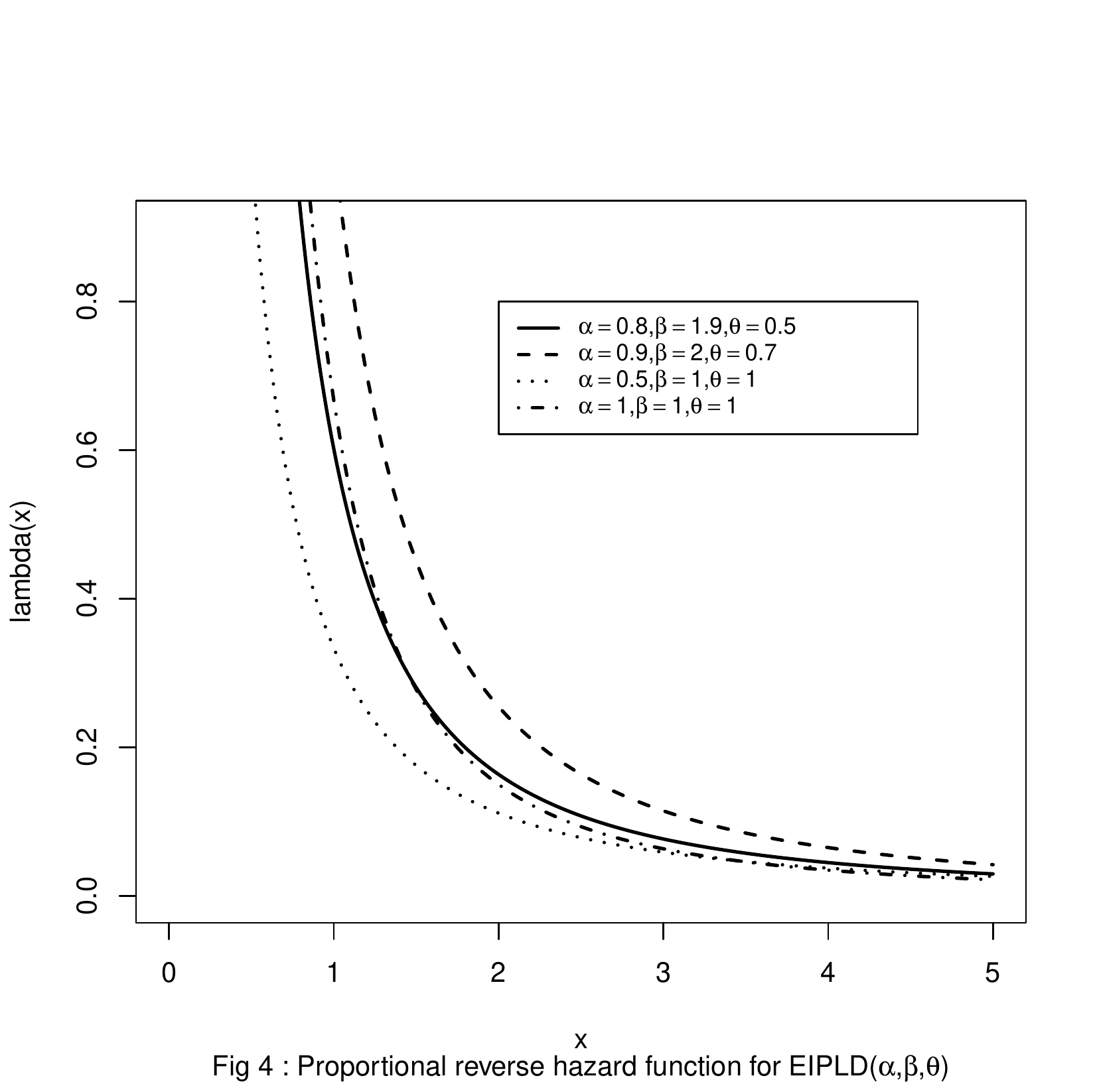}
			
		\end{subfigure}
	\end{figure}

\section{Moments of EIPLD}
Moments play a very important role in probability distributions. Moments are the constants of a population and these constants help in deciding the characteristics of the population and on the basis of these characteristics a population is discussed.\\
\textbf{Theorem 1}: \; If $Z$ be a random variable following EIPL distribution given as in $(4)$, then the $r^{th}$ moment about origin is given by :\\	
\begin{equation}
  		\mu^{'}_{r}=E(Z^r)=(\beta\theta)^{\frac{r}{\alpha}}\sum_{i=0}^{\infty}{\theta-1 \choose i}\frac{1}{(\theta(\beta+1))^{i+1}}\left(i+1-\frac{r}{\alpha}+\theta\beta\right)\Gamma\left({i+1-\frac{r}{\alpha}}\right) \\    \nonumber
\end{equation} 
\textbf{Proof}\\The $r^{th}$ moment about origin of exponentiated inverse power lindley distribution is given by
\begin{eqnarray}
E(Z^{r})&=&\int_{0}^{\infty}z^{r}g_{\theta}(z)dz \nonumber \\
&=&\frac{\alpha\beta^{2}\theta}{1+\beta}\int_{0}^{\infty}z^{r}\left(\frac{1+z^{a}}{z^{2\alpha+1}}\right)e^{\frac{-\beta\theta}{z^{\alpha}}}\left[1+\frac{\beta}{1+\beta}\frac{1}{z^{\alpha}}\right]^{\theta-1}dz 
\end{eqnarray}
Using the binomial series expansion of $\left[1+\frac{\beta}{1+\beta}\frac{1}{z^{\alpha}}\right]^{\theta-1}$ given by
\begin{equation}
\left[1+\frac{\beta}{1+\beta}\frac{1}{z^{\alpha}}\right]^{\theta-1}=\sum_{i=0}^{\infty}{\theta-1\choose i}\left(\frac{\beta}{1+\beta}\frac{1}{z^{\alpha}}\right)^{i} \nonumber
\end{equation}
Equation $(8)$  becomes
\begin{equation}
E(Z^{r})=\frac{\alpha\beta ^2\theta}{1+\beta}\sum_{i=0}^{\infty}{\theta -1\choose i}\left(\frac{\beta}{1+\beta}\right)^i\left[\int_{0}^{\infty}\frac{1}{z^{\alpha(2+i)-r+1}}e^{\frac{-\beta\theta}{z^{\alpha}}}dz    
+\int_{0}^{\infty}\frac{1}{z^{\alpha(i+1)-r+1}}e^{\frac{-\beta\theta}{z^{\alpha}}}dz\right] 
\end{equation}
Let $t=z^{\alpha}$ and using the definition of inverse gamma distribution $(9)$ reduces to
\begin{eqnarray}
E(Z^r)&=&\frac{\beta ^2\theta}{1+\beta}\sum_{i=0}^{\infty}{\theta -1\choose i}\left(\frac{\beta}{1+\beta}\right)^i\left[\frac{\Gamma(2+i-\frac{r}{\alpha})}{(\theta\beta)^{2+i-\frac{r}{\alpha}}}+\frac{\Gamma(i+1-\frac{r}{\alpha})}{(\theta\beta)^{i+1-\frac{r}{\alpha}}}\right]     \nonumber  \\
\mu^{'}_{r}=E(Z^r)&=&(\beta\theta)^{\frac{r}{\alpha}}\sum_{i=0}^{\infty}{\theta-1 \choose i}\frac{1}{(\theta(\beta+1))^{i+1}}\left(i+1-\frac{r}{\alpha}+\theta\beta\right)\Gamma\left({i+1-\frac{r}{\alpha}}\right) 
\end{eqnarray} \\
For $r^{th}$ moment to exist, the contraint $\alpha>r$ must be satisfied.\\
Note that for $\theta=1$, $(10)$ reduces to $r^{th}$ moment of inverse power lindley distribution.
\section{Moment Generating Function}
\textbf{Theorem 2}: \; Let $Z$ be a random variable of EIPLD having pdf $(4)$,then the moment generating function $M_{Z}(t)$ is given by:\\
\begin{equation}
M_{Z}(t)=\sum_{i=0}^{\infty}\sum_{n=0}^{\infty}\frac{t^{n}}{n!}(\beta\theta)^{\frac{n}{\alpha}}{\theta-1 \choose i}\frac{1}{(\theta(\beta+1))^{i+1}}\left(i+1-\frac{n}{\alpha}+\theta\beta\right)\Gamma\left({i+1-\frac{n}{\alpha}}\right) \\ \nonumber
\end{equation}

\textbf{Proof}
\begin{eqnarray}
M_{Z}(t)&=& \int_{0}^{\infty}e^{tz}g_{\theta}(z)dz \nonumber \\
M_{Z}(t)&=& \int_{0}^{\infty}e^{tz}\frac{\alpha\beta^{2}\theta}{1+\beta}\left(\frac{1+z^\alpha}{z^{2\alpha+1}}\right)e^{\frac{-\theta\beta}{z^\alpha}}\left[1+\frac{\beta}{1+\beta}\frac{1}{z^{\alpha}}\right]^{\theta-1}dz 
\end{eqnarray}
Using the binomial series expansion of$\left[1+\frac{\beta}{1+\beta}\frac{1}{z^{\alpha}}\right]^{\theta-1}$ given by:
\begin{equation}
\left[1+\frac{\beta}{1+\beta}\frac{1}{z^{\alpha}}\right]^{\theta-1}= \sum_{i=0}^{\infty}{\theta-1\choose i}\left(\frac{\beta}{1+\beta}\frac{1}{z^{\alpha}}\right)^{i}  \nonumber \\
\end{equation}
Equation $(11)$ becomes	
\begin{eqnarray}
M_{Z}(t)&=&\frac{\alpha\beta^{2}\theta}{1+\beta}\sum_{i=0}^{\infty}{\theta-1 \choose i}\left(\frac{\beta}{1+\beta}\right)^{i}\int_{0}^{\infty}e^{tz}\left( \frac{e^{-\frac{\theta\beta}{z^{\alpha}}}}{z^{\alpha(i+2)+1}} +\frac{e^{-\frac{\theta\beta}{z^{\alpha}}}}{z^{\alpha(i+1)+1}}\right)dz 
\end{eqnarray}
Using $e^{tz}=\sum_{n=0}^{\infty}\frac{t^{n}z^{n}}{n!}$,the equation $(12)$ reduces to
\begin{equation}
M_{Z}(t)=\frac{\alpha\beta^{2}\theta}{1+\beta}\sum_{i=0}^{\infty}\sum_{n=0}^{\infty}{\theta-1 \choose i}\left(\frac{\beta}{1+\beta}\right)^{i}\frac{t^{n}}{n!}\int_{0}^{\infty}\left( \frac{e^{-\frac{\theta\beta}{z^{\alpha}}}}{z^{\alpha(i+2)-n+1}} +\frac{e^{-\frac{\theta\beta}{z^{\alpha}}}}{z^{\alpha(i+1)-n+1}}\right)dz \nonumber
\end{equation}
Let $v=z^{\alpha}$ and using the definition of inverse gamma distribution the above equation becomes
\begin{equation}
M_{Z}(t)=\frac{\alpha\beta^{2}\theta}{1+\beta}\sum_{i=0}^{\infty}\sum_{n=0}^{\infty}{\theta-1 \choose i}\left(\frac{\beta}{1+\beta}\right)^{i}\frac{t^{n}}{n!}\left(\frac{\Gamma(2+i-\frac{n}{\alpha})}{(\theta\beta)^{i+2-\frac{n}{\alpha}}}+\frac{\Gamma(1+i-\frac{n}{\alpha})}{(\theta\beta)^{i+1-\frac{n}{\alpha}}}\right) \nonumber
\end{equation}
Thus the moment generating function of EIPLD is given by:
\begin{equation}
M_{Z}(t)=\sum_{i=0}^{\infty}\sum_{n=0}^{\infty}\frac{t^{n}}{n!}(\beta\theta)^{\frac{n}{\alpha}}{\theta-1 \choose i}\frac{1}{(\theta(\beta+1))^{i+1}}\left(i+1-\frac{n}{\alpha}+\theta\beta\right)\Gamma\left({i+1-\frac{n}{\alpha}}\right)  \nonumber
\end{equation}

\section{Quantile Function}
\textbf{Theorem 3}: \; If Z follows EIPLD$(\alpha,\beta,\theta)$ then the quantile function of Z is 
\begin{equation}
Q(u)=\left[-1-\frac{1}{\beta}-\frac{1}{\beta}W_{-1}\left(-u^{\frac{1}{\theta}}(1+\beta)e^{-(1+\beta)}\right)\right]^{-\frac{1}{\alpha}} \nonumber
\end{equation}
where $u\in(0,1)$ and $W_{-1}$ denote the negative branch of Lambert $W$ function.\\

\textbf{Proof} : \;
The quantile function denoted by $Q(u)$ is the root of the equation
\begin{equation}
\left[\left(1+\frac{\beta}{1+\beta}\frac{1}{Q(u)^\alpha}\right)e^\frac{-\beta}{Q(u)^\alpha}\right]^\theta=u \; \;  ;0<u<1  
\end{equation}
Multiplying $(13)$ both sides by $e^{-1-\beta}$ we get,
\begin{equation}
-\left[1+\beta+\frac{\beta}{Q(u)^{\alpha}}\right]e^{-\left(1+\beta+\frac{\beta}{Q(u)^{\alpha}}\right)}=-(1+\beta)u^{\frac{1}{\theta}}e^{-(1+\beta)}    \nonumber
\end{equation}
Using the Lambert W function which is the solution of the equation $W(z)e^{W(z)}=z$,where z is a complex number ,we have
\begin{equation}
W\left(-u^{\frac{1}{\theta}}e^{-(1+\beta)}(1+\beta)\right)=-\left(1+\beta+\frac{\beta}{Q(u)^{\alpha}}\right)     \nonumber
\end{equation}
The negative Lambert $W$ function of the real argument $-u(1+\beta)e^{1+\beta}$ is
\begin{equation}
W_{-1}\left(-u^{\frac{1}{\theta}}e^{-(1+\beta)}(1+\beta)\right)=-\left(1+\beta+\frac{\beta}{Q(u)^{\alpha}}\right)     \nonumber
\end{equation}
which upon solving for $Q(u)$ results in
\begin{equation}
Q(u)=\left[-1-\frac{1}{\beta}-\frac{1}{\beta}W_{-1}\left(-u^{\frac{1}{\theta}}(1+\beta)e^{-(1+\beta)}\right)\right]^{-\frac{1}{\alpha}}     \nonumber
\end{equation}
Using above equation the quartiles of the exponentiated inverse power lindley distribution can be determined.Median of exponentiated inverse power lindley distribution is given by
\begin{equation}
Q\left(\frac{1}{2}\right)=\left[-1-\frac{1}{\beta}-\frac{1}{\beta}W_{-1}\left(-\left(\frac{1}{2}\right)^{\frac{1}{\theta}}(1+\beta)e^{-(1+\beta)}\right)\right]^{-\frac{1}{\alpha}}     \nonumber
\end{equation}

\section{Renyi entropy}
Entropies quantify the diversity, uncertainty, or randomness of a system. The Renyi entropy is named after Alfred Renyi in the context of fractal dimension estimation, the Renyi entropy forms the basis of the concept of generalized dimensions. For a given probability distribution, Renyi entropy is given by:\\
 
\begin{equation}
e(\gamma)=\frac{1}{1-\gamma}log\left[\int{g^{\gamma}(z)dz}\right]     \nonumber
\end{equation} 
where $\gamma>0 \; and \; \gamma\neq1$

\begin{equation}
e(\gamma)=\frac{1}{1-\gamma}log\left[\int_{0}^{\infty}\frac{\alpha\beta^{2}\theta}{1+\beta}\left(\frac{1+z^{\alpha}}{z^{2\alpha+1}}\right)e^{-\frac{\beta\theta}{z^{\alpha}}}\left(1+\frac{\beta}{1+\beta}\frac{1}{z^{\alpha}}\right)^{\theta-1}dz\right]^{\gamma}    \nonumber
\end{equation}
Using the binomial series expansion $(1+z)^{n}=\sum_{i=0}^{\infty}{n\choose i}z^{i}$ for $\left[\left(1+\frac{\beta}{1+\beta}\frac{1}{z^{\alpha}}\right)^{\theta-1}\right]^{\gamma}$  and $\left(1+\frac{1}{z^{\alpha}}\right)^{\gamma}$ \\ 
 we have,
\begin{equation}
e(\gamma)=\frac{1}{1-\gamma}log\left(\frac{\alpha\beta^{2}\gamma}{1+\beta}\right)^{\gamma}\sum_{i=0}^{\infty}{\gamma\theta-\gamma \choose i}\left(\frac{\beta}{1+\beta}\right)^{i}\sum_{j=0}^{\infty}{\gamma\choose j}\int_{0}^{\infty}e^{-\frac{\gamma\beta\alpha}{z^{\alpha}}}\frac{1}{z^{\alpha(i+j+\gamma)+\gamma}}dz   
\end{equation}
\;

Letting $z^{\alpha}=t$ and using the definition of inverse gamma distribution $(14)$ reduces to 
\begin{equation}
e(\gamma)=\frac{1}{1-\gamma}log\left( \frac{\alpha^{\gamma-1}\beta^{2\gamma}\theta^{\gamma}}{(1+\beta)^{\gamma}}\right)\sum_{i=0}^{\infty}\sum_{j=0}^{\infty}{\gamma\theta-\gamma \choose i}{\gamma \choose j}\left(\frac{\beta}{1+\beta}\right)^{i} \frac{\Gamma\left(i+j+\gamma+\frac{\gamma-1}{\alpha}\right)}{(\theta\beta\gamma)^{i+j+\gamma+\frac{\gamma-1}{\alpha}}}    \nonumber
\end{equation}

\section{Distribution of Order Statistic}

The pdf of $k^{th}$ order statistic is given by:

\begin{equation}
g(z)_{k}=\frac{n!g_{\theta}(z)}{(k-1)!(n-k)!}(G(z))^{k-1}(1-G(z))^{n-k}     \nonumber
\end{equation}
where $g_{\theta}(z)$ and $G(z)$ denotes the pdf and cdf respectively. \\
Using the binomial series expansion $(1-z)^{n}=\sum_{i=0}^{\infty}{n \choose i}(-1)^{i}(z)^{i}$, we get
\begin{eqnarray}
g(z)_{k}&=&\frac{n!g_{\theta}(z)}{(k-1)!(n-k)!}\sum_{i=0}^{\infty}{n-k \choose i}(-1)^{i}(G(z))^{i+k-1}g_{\theta}(z)  \nonumber  \\
&=&\frac{n!}{(k-1)!(n-k)!}\sum_{i=0}^{\infty}(-1)^{i}{n-k \choose i}\frac{\alpha\beta^{2}\theta}{1+\beta}\left(\frac{1+z^{\alpha}}{z^{2\alpha+1}}\right)e^{-\frac{\theta\beta}{z^{\alpha}}}\left(1+\frac{\beta}{1+\beta}\frac{1}{z^{\alpha}}\right)^{\theta i+\theta k-1}   \nonumber \\
&=&\frac{\alpha\beta^{2}\theta}{1+\beta}\frac{n!}{(k-1)!(n-k)!}\sum_{i=0}^{\infty}\sum_{j=0}^{\infty}(-1)^{i}{n-k \choose i}{\theta i+\theta k-1 \choose j} \nonumber\\ &&\left(\frac{\beta}{1+\beta}\right)^{j}\left(\frac{1}{z^{\alpha}}\right)^{j}\left(\frac{1+z^{\alpha}}{z^{2\alpha+1}}\right)e^{-\frac{\theta\beta}{z^{\alpha}}}      \nonumber  \\
&=&\frac{\alpha\theta\beta^{j+2}n!}{(1+\beta)^{j+1}(k-1)!(n-k)!}\sum_{i=0}^{\infty}\sum_{j=0}^{\infty}(-1)^{i}{n-k \choose i}{\theta i+\theta k-1 \choose j}  \nonumber \\
&&\left(\frac{1+z^{\alpha}}{z^{(j+2)\alpha+1}}\right)e^{-\frac{\theta\beta}{z^{\alpha}}}   \nonumber
\end{eqnarray}

\section{Stochastic Ordering}

Stochastic ordering of a continous random variable is an important tool for judging the comparative behavior. A random variable $Z$ is said to be greater than $Y$ in the:\\
\;
(a) Stochastic order $(Y\leq_{st}Z)$ if $G{_{Z}}(z)\leq G_{_{Y}}(z) \forall z$\\
(b) Hazard rate order $(Y\leq_{hr}Z)$ if $h_{Z}(z)\leq h_{Y}(z) \forall z$\\
(c) Mean residual order $(Y\leq_{mlr}Z)$ if $m_{Z}(z)\leq m_{Y}(z)\forall z$\\
(d) Likelihood ratio order $(Y\leq_{lr}Z)$ if $\frac{g_{Z}(z)}{g_{Y}(z)}$ is an increasing function of $z$.\\

\;
The following results(Shaked and Shanthikumar 1994) are well known:\\
\begin{center}
	{$Z\leq_{lr}Y \Rightarrow Z\leq_{hr}Y \Rightarrow Z\leq_{mrl}Y $}\\
	{$\Downarrow$}               \\
	{$Z\leq_{st}Y$}
\end{center}

The following theorem shows that EIPLD is ordered with respect to "likelihood ratio " ordering.\\

\;
\textbf{Theorem 4} :
Let $Y \sim EIPLD (\alpha_{1},\beta_{1},\theta_{1})$ and $Z \sim EIPLD(\alpha_{2},\beta_{2},\theta_{2})$. .If $\beta_{1}=\beta_{2}$ and $\theta_{2}\geq\theta_{1}$  (or if $\beta_{2}\geq\beta_{1}$ and $\theta_{1}=\theta_{2}$ ) then $(Y\leq_{lr}Z) \;  \forall \; z$.

\textbf{Proof}
We have \\

\begin{eqnarray}
\frac{g_{Z}(z)}{g_{Y}(z)}&=&\frac{\alpha_{2}\beta_{2}^{2}\theta_{2}(1+\beta_{1})}{\alpha_{1}\beta_{1}^{2}\theta_{1}(1+\beta_{2})}\left(\frac{1+z^{\alpha_{2}}}{1+z^{\alpha_{1}}}\right)\left(\frac{z^{2\alpha_{1}+1}}{z^{\alpha_{2}+1}}\right) \nonumber \\
&&exp{-\left(\frac{\theta_{2}\beta_{2}}{z^{\alpha_{2}}}-\frac{\theta_{1}\beta_{1}}{z^{\alpha_{1}}}\right)}  \left[\frac{\left(1+\frac{\beta_{2}}{1+\beta_{2}}\frac{1}{z^{\alpha_{2}}}\right)^{\theta_{2}-1}} {\left(1+\frac{\beta_{1}}{1+\beta_{1}}\frac{1}{z^{\alpha_{1}}}\right)^{\theta_{1}-1}}\right]  \nonumber 
\end{eqnarray} 
Setting $\alpha_{1}=\alpha_{2}$\\
\textbf{Case 1}: for $\beta_{1}=\beta_{2}$ and $\theta_{2}\geq\theta_{1}$ 
we obtained $\frac{d}{dz}\left(\frac{g_{Z}(z)}{g_{Y}(z)}\right)$ as an increasing function of $z$. \\
\textbf{Case 2}: $\beta_{2}\geq\beta_{1}$ and $\theta_{1}=\theta_{2}$ we obtained $\frac{d}{dz}\left(\frac{g_{Z}(z)}{g_{Y}(z)}\right)$ as an increasing function of $z$.\\
This implies that $Y\leq_{lr}Z \;  \forall \; z$.Hence $Y\leq_{hr}Z$ , $Y\leq_{mrl}Z$ and $Y\leq_{st}Z$. \\

\section{Maximum Likelihood Estimation}
Let $z_{1},.....z_{n}$ be a random sample of size $n$ from the EIPLD. The log-likelihood equation is given by\\
\begin{eqnarray}
L(\alpha,\beta,\theta |z) &=& n[log\alpha+2log\beta+log\theta-log(1+\beta)]+\sum_{i=1}^{n}log(1+z_{i}^{\alpha})-(2\alpha+1)\sum_{i=1}^{n}log(z_{i}) \nonumber \\
&-&\theta\beta\sum_{i=1}^{n}z_{i}^{-\alpha}  (\theta-1)\sum_{i=1}^{n}log\left[1+\frac{\beta}{1+\beta}\frac{1}{z_{i}^{\alpha}}\right] \nonumber
\end{eqnarray}
The maximum likelihood estimates $\hat{\alpha},\hat{\beta},\hat{\theta}$ for $\alpha,\beta,\theta$ are :
\begin{eqnarray}
\frac{\partial}{\partial\alpha}L(\alpha,\beta,\theta |z)&=&\frac{n}{\alpha}+\sum_{i=1}^{n}\frac{z_{i}^{\alpha}logz_{i}}{1+z_{i}^{\alpha}}-2\sum_{i=1}^{n}log(z_{i})+\theta\beta\sum_{i=1}^{n}(z_{i}^{-\alpha})logz_{i} \nonumber \\
 &&-(\theta-1)\sum_{i=1}^{n}\left[\frac{\frac{\beta}{1+\beta}(z_{i}^{-\alpha})log(z_{i})}{1+\frac{\beta}{1+\beta}\frac{1}{z_{i}}^{\alpha}}\right]=0 \nonumber\\
\frac{\partial}{\partial\beta}L(\alpha,\beta,\theta |z)&=&\frac{n(2+\beta)}{\beta(1+\beta)}-\theta\sum_{i=1}^{n}z_{i}^{-\alpha} \nonumber \\
  &&+(\theta-1)\sum_{i=1}^{n}\left[\frac{1}{1+\frac{\beta}{1+\beta}\frac{1}{z_{i}^{\alpha}}}\right] \left[\frac{1}{z_{i}^{\alpha}}\frac{1}{(1+\beta^{2})}\right]=0 \nonumber\\
\frac{\partial}{\partial\theta}L(\alpha,\beta,\theta |z)&=&\frac{n}{\theta}-\beta\sum_{i=1}^{n}z_{i}^{-\alpha}+\sum_{i=1}^{n}log\left[1+\frac{\beta}{1+\beta}\frac{1}{z_{i}^{\alpha}}\right]=0 \nonumber
\end{eqnarray}
The above non-linear system of equations are solved by numerical iteration technique and maximum likelihood estimates are obtained.\\
Since the maximum likelihood estimates for $\alpha,\beta,\theta$ are not in closed form we use the large sample behaviour of maximum likelihood estimators to obtain the confidence intervals for model parameters. The asymptotic sampling distribution of $(\hat{\alpha},\hat{\beta},\hat{\theta})$ is $N[(\alpha,\beta,\theta),\Delta^{-1}]$ where $\Delta$ is the observed Fisher information matrix given by:

\[
\Delta=
\begin{Bmatrix}
\frac{\partial^{2}L}{\partial\alpha^{2}} & \frac{\partial^{2}L}{\partial\alpha\partial\beta}  & \frac{\partial^{2}L}{\partial\alpha\partial\theta}\\
\frac{\partial^{2}L}{\partial\beta\partial\alpha}&\frac{\partial^{2}L}{\partial\beta^{2}}&\frac{\partial^{2}L}{\partial\beta\partial\theta}\\                    
\frac{\partial^{2}L}{\partial\theta\partial\alpha}&\frac{\partial^{2}L}{\partial\theta\partial\beta} & \frac{\partial^{2}L}{\partial\theta^{2}}
\end{Bmatrix}
\] or
\[
\Delta=
\begin{Bmatrix}
\widehat{Var(\alpha)} & \widehat{Cov(\alpha,\beta)} & \widehat{Cov(\alpha.\theta)}\\
\widehat{Cov(\beta,\alpha)} & \widehat{Var(\beta)} & \widehat{Cov(\beta,\theta)}\\                    
\widehat{Cov(\theta,\alpha)} & \widehat{Cov(\theta,
\beta)} &\widehat{ Var(\theta)}
\end{Bmatrix}
\]  

The second order derivatives for the parameters of EIPLD$(\alpha,\beta,\theta)$ exist and are derived as:
\begin{eqnarray}
\frac{\partial^{2}L}{\partial\alpha^{2}}&=&-\frac{n}{\alpha^{2}}+\sum_{i=1}^{n}\frac{z_{i}^{\alpha}(logz_{i})^{2}}{(1+z_{i}^{\alpha})^{2}}-\theta\beta\sum_{i=1}^{n}z_{i}^{-\alpha}logz_{i} \nonumber \\ &&+\;  (\theta-1)\left(\frac{\beta}{1+\beta}\right) \; \sum_{i=1}^{n}\frac{z_{i}^{-\alpha}(logz_{i})^{2}}{\left(1+\frac{\beta}{1+\beta}\frac{1}{z_{i}^{\alpha}}\right)^{2}} \nonumber \\
\frac{\partial^{2}L}{\partial\beta^{2}}&=&\frac{-2n}{\beta^{2}}+\frac{n}{(1+\beta)^{2}} \nonumber \\
&&- \; \frac{(\theta-1)}{(1+\beta)^{4}} \; \sum_{i=1}^{n}\left[\left( \frac{1}{1+\frac{\beta}{1+\beta}\frac{1}{z_{i}^{\alpha}}}\right)\frac{1}{z_{i}^{\alpha}}\right]^{2}\left( 3+2\beta(1+z_{i})\right) \nonumber \\
\frac{\partial^{2}L}{\partial\theta^{2}}&=& \frac{-n}{\theta^{2}} \nonumber \\
\frac{\partial^{2}L}{\partial\alpha\partial\theta}&=&\beta \; \sum_{i=1}^{n}(z_{i}^{-\alpha})logz_{i}-\frac{\beta}{1+\beta} \; \sum_{i=1}^{n} \left[\frac{z_{i}^{-\alpha}logz_{i}}{1+\frac{\beta}{1+\beta}\frac{1}{z_{i}^{\alpha}}} \right] \nonumber \\
\frac{\partial^{2}L}{\partial\beta\partial\theta}&=& -\frac{\beta}{(1+\beta)^{3}} \; \sum_{i=1}^{n}\left(\frac{1}{z_{i}^{\alpha}}\right)^{2} \nonumber \\ \frac{\partial^{2}L}{\partial\beta\partial\alpha}&=&\theta \sum_{i=1}^{n}z_{i}^{-\alpha}logz_{i} - \frac{(\theta-1)}{(1+\beta)^{2}} \, \sum_{i=1}^{n} \frac{z_{i}^{-\alpha}logz_{i}}{\left[1+\frac{\beta}{1+\beta}\frac{1}{z_{i}^{\alpha}}\right]^{2}} \nonumber 
\end{eqnarray}
The solutions of the above equations yield the asymptotic variance covariance of ML estimators for $\alpha,\beta,\theta$.
The asymptotic confidence intervals for $\alpha,\beta,\theta$ are determined as:

$\hat{\alpha}+Z_{\frac{\alpha}{2}}\sqrt{\hat{Var(\alpha)}} \;  , \;\; \; \hat{\beta}+Z_{\frac{\alpha}{2}}\sqrt{\hat{Var(\beta)}} \; , \; \; \; \hat{\theta}+Z_{\frac{\alpha}{2}}\sqrt{\hat{Var(\theta)}}$

\section{Data analysis}
In this section, real data analysis is performed to illustrate the applicability of Exponentiated Inverse Power Lindley distribution (EIPLD). The data set given in Table 1 represents the active repair times (hr) for an airborne communication transceiver.This data has been widely used by various authors and were initially used by Jorjensen(1982).

	\begin{table}
	\caption{Active repair times}
	\begin{center}

		\begin{tabular}{|c|c|c|c|c|c|c|c|}

			\hline
			0.50 & 0.60 & 0.60 &0.70 & 0.70 &0.70 &0.80 &0.80  \\ [1ex] \hline
			1.00 & 1.00 & 1.00 &1.00 &1.10& 1.30 & 1.50 &1.50 \\ [1ex] \hline
			1.50 & 1.50 &2.00 &2.00 & 2.20 &2.50 & 2.70 &3.00 \\ [1ex] \hline
			3.00 & 3.30 & 4.00 & 4.00 & 4.50 & 4.70 &5.00 &5.40 \\ [1ex] \hline
			5.40 & 7.00 &7.50 & 8.80 & 9.00 & 10.20 & 22.00 & 24.50\\ [1ex] \hline
		
		\end{tabular}
	\end{center}
   \end{table}

\newpage
The applicability of EIPLD is demonstrated by using some statistical tools suchs as Kolmogrov-Smirnov statistic, Akaike information criterion(AIC) defined by $-2logL+2q$, Bayesian information criterion(BIC) defined by $-2logL+qlog(n)$, where q is the number of estimated parameters and n is the sample size, and are compared with other distributions. AIC and BIC values estimates the quality of each model relative to each of the other models.The MLEs of the parameters are given in Table 3 and the statistical values mentioned above are computed and are given in Table 2. These values indicate that the proposed distribution fits well to the data compared to other tested distributions.
     
	\begin{table}[h!]
 	\caption{Comparison criterion}
 \begin{center}
 	\begin{tabular}{|c|c|c|c|c|c|}
 		\hline
 		Model    & Model      &   AIC    &   BIC    &  -Logl  & K-S statistic \\
 		[1ex]\hline	
 		Exponentiated Inverse Power Lindley Distribution & EIPLD &  184.91  &  189.97  &  89.45  &    0.09537    \\[1ex] \hline
 		Exponentiated Power Lindley Distribution & EPLD  & 186.5721 & 191.6387 & 90.2861 &     0.909     \\[1ex] \hline
 		Power Lindley Distribution &PLD & 195.8854 & 199.2631 & 95.9427 &    0.1596     \\[1ex]  \hline
 		Generalised Lindley Distribution & GLD & 199.8218 & 203.1995 & 97.9109 &    0.1410     \\[1ex] \hline
 		Lindley Distribution &LD & 199.8218 & 203.1995 & 97.9109 &    0.1907     \\[1ex] \hline
 		Exponentiated Exponential & EE  & 194.9158 & 198.2936 & 95.4579 &    0.1334     \\[1ex] \hline
 		Weibull Distribution &WD  & 195.0227 & 198.4005 & 95.5114 &    0.1540     \\[1ex]  \hline      
 	\end{tabular}	
 \end{center}
\end{table}
   
   \begin{table}
	\caption{Parameters MLES}
\begin{center}
	\begin{tabular}{|c|c|c|c|}
		\hline
		Model & $\alpha$ &$\beta$ &$\theta$ \\[1ex]\hline
		EIPLD & 1.20167 & 25.94112 & 0.06205 \\[1ex]
		EPLD  &30.8299 &3.5472 &0.2901 \\[1ex]
		PLD &-&0.5867 &0.7988  \\[1ex]
		GLD &0.7460 &0.3588 &- \\[1ex]
		EE  &1.1137 &0.2678 &- \\[1ex]
		WD  &- &0.2688 &0.9604  \\[1ex]\hline
	\end{tabular}
\end{center}
\end{table}

\newpage
\section{Simulation Study}
In this section, we study the performance of ML estimators for different sample sizes $(n=25,50,100,200,300,500)$. We have employed the inverse CDF technique for data simulation for EIPLD using R software. The process was repeated 500 times for calculation of bias, variance and MSE.  For two parameter combinations of EIPLD, decreasing trend is being observed in average bias, variance and MSE as we increase the sample size. Hence, the performance of ML estimators is quite well, consistent in case of EIPLD. 
\begin{figure}[h!]
	\begin{center}
		\includegraphics[scale=1.0]{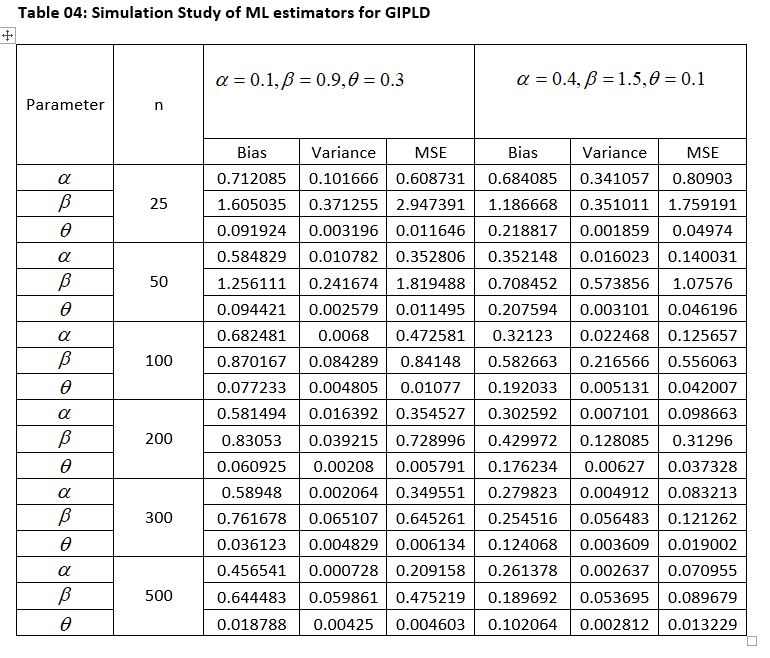}
	\end{center}
\end{figure}
\newpage
\section{Conclusion}
In this article, a three parameter distribution called as exponentiated inverse power lindley distribution has been proposed. Some statistical properties such as moments,moment generating function,quantile,stochastic ordering, renyi entropy of the proposed distribution has been discussed. The parameter estimation is approached by method of maximum likelihood estimation. Confidence intervals for the model parameters are also derived. An application of EIPLD to the real data set shows that this model could provide a better fit than other compared models. Simulation study was performed using different values of parameters and MSE,bias has been obtained.

\paragraph{References}
\begin{itemize}
\item Alkarni, S. 2015. Extended Power Lindley distribution-A new Statistical model for non-monotone survival data. European journal of statistics and probability 3(3):19-34.
\item Ashour, s., and Eltehiwy, M. A. 2015. Exponentiated Power Lindley distribution. J Adv Res 6(6): 895-905.
\item Bakouch, H.S., Al-Zaharani, B., Al-Shomrani, A., Marchi, V., and Louzad, F. 2012. An  extended Lindley distribution. Journal of the Korean Statistical Society 41(1):75-85.
\item Barco, Mazucheli and Janerio. 2016. The inverse power lindley distribution. Communication in statistics-Simulation and Computation. 46(8):6308-6323.
\item Ghitany, M., Al-Mutairi, D., Balakrishnan, N., and Al-Enezi, I. 2013. Power Lindley distribution and associated inference. Computational Statistics and Data Analysis 64:20-33.
\item Ghitany, M., Atieh, B., Nadarajah,S. 2008. Lindley distribution and its Application. Mathematics Computing and Simulation 78(4): 493-506.
\item Jorgensen, B. 1982. Statistical properties of the generalized inverse Gaussian distribution. New York: Springer-Verlag.
\item Lindley, D. V. 1958. Fiducial distributions and Bayes’ theorem. Journal of the Royal  Statistical Society Series B 20(1): 102-107.
\item Mudholkar, G., and Srivastava, D. 1993. Exponentiated weibull family for analyzing bathtub failure rate data. IEEE Transactions on Reliability 42, 299-302.
\item Nadarajah, S., Bakouch, H. S., and Tahmasbi, R. 2011. A generalized Lindley distribution. Sankhya B 73(2): 331-359.
\item Shaked, M., and Shanthikumar, J. G. 1994. Stochastic orders and their applications. New York: Academic Press.
\item Shanker, R., and Mishra, A. 2013. A two-parameter Lindley distribution. Statistics in Transition-new series 14 (1): 45-56.

\item Sharma, V., Singh, S., Singh, U., and Agiwal, V. 2015. The inverse Lindley distribution- A stress-strength reliability model with applications to head and neck cancer data. Journal  of Industrial and Production Engineering 32(3): 162-173.
\item Sharma, V., Singh, S., Singh, U., and Merovci, F.  2015.  The generalized inverse Lindley distribution: A new inverse statistical model for the study of upside-down bathtub survival data. Commun Stat Theory Methods, preprint.
\item Xie, M., Goh, T., Tang, Y. 2002. A modified weibull extension with bathtub-shaped failure rate function. Reliability Engineering and Systems Safety 76, 279-285.
\item Xie, M., and Lai, C. 1996. Reliability analysis using an additive weibull model with bathtub shaped failure rate function. Reliability Engineering and Systems Safety 52, 87-93.
\item Zakerzadeh, H., and Dolati, A .2009. Generalized Lindley distribution. Journal of Mathematical extension 3(2): 13-25.
\end{itemize}

\end{document}